\documentclass[a4paper,reqno]{article}

\usepackage{amssymb}

\def\bdn{\begin{description}}
\def\edn{\end{description}}
\def\biz{\begin{itemize}}
\def\eiz{\end{itemize}}   
\def\ben{\begin{enumerate}}
\def\een{\end{enumerate}}
\def\bt{\begin{center}\begin{tabular}}
\def\et{\end{tabular}\end{center}}
\def\be{\begin{equation}}
\def\ee{\end{equation}}
\def\bar{\begin{eqnarray}}
\def\ear{\end{eqnarray}}

\def\o{\`o}
\def\u{\`u}

\def\avalori{\longrightarrow}

\def\de{\partial}

\def\ov#1{\overline{#1}}

\def\unione{\cup}
\def\intersez{\cap}

\newcommand{\bdm}{\begin{displaymath}}
\newcommand{\edm}{\end{displaymath}}
\newcommand{\p}{\partial}

\def\invstackrel#1#2{\mathord{\mathop{#1}\limits_{#2}}}

\newsavebox{\junk}
\newsavebox{\junket}

\def\bddots{\mathinner{\mkern1mu\raise1pt\hbox{.}\mkern2mu
        \raise4pt\hbox{.}\mkern2mu\raise7pt\vbox{\kern7pt\hbox{.}}\mkern1mu}}

\def\bkB{{\rm I\kern-.17em B}}
\def\bkC{{\rm \kern.24em
            \vrule width.05em height1.4ex depth-.05ex
            \kern-.26em C}}
\def\bkD{{\rm I\kern-.17em D}}
\def\bkE{{\rm I\kern-.17em E}}
\def\bkF{{\rm I\kern-.17em F}}
\def\bkG{{\rm \kern.24em
            \vrule width.05em height1.4ex depth-.05ex
            \kern-.26em G}}
\def\bkH{{\rm I\kern-.22em H}}
\def\bkI{{\rm I\kern-.22em I}}
\def\bkJ{{\rm \kern.19em
            \vrule width.02em height1.5ex depth0ex
            \kern-.20em J}}
\def\bkK{{\rm I\kern-.22em K}}
\def\bkL{{\rm I\kern-.17em L}}
\def\bkM{{\rm I\kern-.22em M}}
\def\bkN{{\rm I\kern-.20em N}}
\def\bkO{{\rm \kern.24em
            \vrule width.05em height1.4ex depth-.05ex
            \kern-.26em O}}
\def\bkP{{\rm I\kern-.17em P}}
\def\bkQ{{\rm \kern.24em
            \vrule width.05em height1.4ex depth-.05ex
            \kern-.26em Q}}
\def\bkR{{\rm I\kern-.17em R}}
\def\bkT{{\rm \kern.24em
            \vrule width.02em height1.5ex depth 0ex
            \kern-.27em T}}
\def\bkU{{\rm \kern.30em
            \vrule width.02em height1.47ex depth-.05ex
            \kern-.32em U}}
\def\bkZ{{\rm Z\kern-.32em Z}}

\def\R{\bkR}
\def\N{\bkN}

\def\E{\bkE}

\newtheorem{Def}{Definition}[section]
\newtheorem{teo}[Def]{Theorem}

\newtheorem{prop}[Def]{Proposition}

\newtheorem{assumption}[Def]{Assumption}

\newtheorem{lem}[Def]{Lemma}

\newtheorem{oss}[Def]{Remark}

\newenvironment{Dim}[0]{{\bf Proof}} {}
\newenvironment{Dimo}[0]{{\bf{Proof of }}} 

\def\o{\overline}
\def\u{\underline}

\def\*{\star}

\def\t{\tilde}

\def\N{{\mathbb N}}

\def\L1{\Lambda_B^1([0,T])}
\def\L2{\Lambda_B^2([0,T])}
\def\l{\lambda}
\def\L{\Lambda}

\def\M1{M^1_B([0,T])}
\def\M2{M^2_B([0,T])}
\def\K{{\cal K}}
\def\D{{\cal D}}

\def\Y{{\cal Y}}

\def\M{{\cal M}}
\def\N{{\cal N}}

\def\F{{\cal F}}

\def\s{\sigma}
\def\p{\prime}
\def\u{\underline}
\def\o{\overline}
\def\t{\tilde}

\def\var{\varepsilon}

\newcommand{\cvd}{\begin{flushright}
\rule[5pt]{5pt}{5pt}
\end{flushright}}

\date{}
\begin{document}
\title{Lipschitzian Estimates in Discrete-Time Constrained Stochastic Optimal Control}
\maketitle
\begin{center}
\author{Marco Papi\\
\vspace*{5pt}
\begin{scriptsize}
Istituto per le Applicazioni del Calcolo ``M.Picone'',\\
 V.le del Policlinico 137, I-00161 Roma (Italy) \\
and\\ 
Dipartimento di Matematica di Roma ``Tor Vergata'',\\
Via della Ricerca Scientifica, 00133, Roma (Italy),\\
e-mail: papi@iac.rm.cnr.it
\end{scriptsize}} 
\end{center}
\begin{center}
and
\end{center}
\begin{center}
\author{Simone Sbaraglia\\
\vspace*{5pt}
\begin{scriptsize}
Istituto per le Applicazioni del Calcolo ``M.Picone'',\\
 V.le del Policlinico 137, I-00161 Roma (Italy) \\
and\\ 
Dipartimento di Metodi e Modelli Matematici di Roma ``La Sapienza'',\\
Via Scarpa 16, 00161, Roma (Italy),\\
e-mail: sbaragli@iac.rm.cnr.it
\end{scriptsize}} 
\end{center}
\vspace*{2pt}

\begin{abstract}
This paper is devoted to the analysis of a finite horizon discrete-time stochastic optimal control
problem, in presence of constraints.  We study the regularity of the value function which comes from
the dynamic programming algorithm. We derive accurate estimates of the Lipschitz constant of the
value function, by means of a regularity result of the multifunction that defines the admissible
control set.\\
In the last section we discuss an application to an optimal asset-allocation problem.
\end{abstract}

\begin{quote}
\begin{small}
{\bf Key words.} Optimal control, Dynamic Programming, State constraints, Lipschitz regularity,
asset-allocation, Multifunctions.\\
{\bf AMS subject classifications.} 49L20, 49N60, 32A12, 37N40, 37N35.
\end{small}
\end{quote}

\vspace*{20pt}

%%%%%%%%%%%%%%%%%%%%%%%%%%%%%%%%%%%%%%%%%%%%%%%%%%%%%%%%%%%%%%%%%%%%%%%%%%%%%%%%%%%%%%%%%%%
%\Section{Introduction}
%%%%%%%%%%%%%%%%%%%%%%%%%%%%%%%%%%%%%%%%%%%%%%%%%%%%%%%%%%%%%%%%%%%%%%%%%%%%%%%%%%%%%%%%%%%

This paper is devoted to the analysis of a general {\em Finite Horizon Discrete-Time Stochastic
Optimal Control Model}, with inequality constraints.
The aim of this paper is to give a method to estimate the Lipschitz constant of the
value function obtained via the classical dynamic programming algorithm.\\

The regularity of the value function is related to the regularity of the {\em marginal function}
\cite{AUBIN}, \cite{VINTER}, and, as is proved in \cite{AUBIN}, the regularity of the marginal 
function is connected to the regularity of the multifunction which defines the set of admissible
controls. Therefore, our main result concerns the Lipschitz regularity, with respect to the
Hausdorff metric \cite{BARNSLEY}, of this multifunction.\\

Set-valued maps are widely used in optimal control, differential games and their applications to
mathematical economics and finance, see \cite{DEBREU}, and \cite{ROCKAFELLAR-LIP}.
In many cases these multifunctions are defined by means of inequality constraints for a set of
functions defined over a manifold that represents the control space. Unfortunately this
manifold is usually non-regular, as in the financial application presented in Section
\ref{application}. Furthermore, the constraint functions may lose the regularity or be dependent on
each other, at some point.\\
To overcome this difficulties, we allow for Lipschitz manifolds and Lipschitz constraint functions,
provided that the set where either the manifold or the constraints are not regular, or the
constraints are dependent on each other, can be approximated by points where both the manifold and
the constraint function are regular and the constraints are independent on each other.\\

We use a quantitative formulation of the Implicit Function Theorem that provides an
estimate of the neighborhoods where the implicit function is defined.\\
Our study is carried out in the discrete-time case, because of its high computational relevance,
since the regularity properties of the value function can be used to derive a-priori error estimates
and convergence results of numerical schemes.\\

The outline of the paper is as follows. Section \ref{MODEL} introduces the general framework of the
optimization model and recalls the classical discrete-time DP algorithm. Section \ref{HAUSDORFF}
provides the definition of the Hausdorff metric and some related results, Section \ref{notations}
provides basic notations and definitions, Section \ref{LIPSCHITZEST} establishes the main regularity
results about multifunctions and a Lipschitzian estimate for the value function related to the
dynamic programming algorithm. Finally, Section \ref{application} applies these results to an
optimal asset-allocation problem with regulatory constraints.

%%%%%%%%%%%%%%%%%%%%%%%%%%%%%%%%%%%%%%%%%%%%%%%%%%%%%%%%%%%%%%%%%%%%%%%%%%%%%%%%%%%%%%%%%%%
\section{The Dynamic Programming Algorithm}\label{MODEL}
%%%%%%%%%%%%%%%%%%%%%%%%%%%%%%%%%%%%%%%%%%%%%%%%%%%%%%%%%%%%%%%%%%%%%%%%%%%%%%%%%%%%%%%%%%%

In this section we present the model which is the subject of our study in next
sections. We consider the following discrete-time controlled dynamical system:
\be
\label{DS}
\left\{
\begin{array}{l}
x_{k+1}=f_k(x_k,u_k,y_k),\;\;\;\;k=0,\ldots,N-1\qquad\qquad\\
\\
x_0=x\in X_0,\qquad \qquad \qquad \qquad \qquad\qquad \qquad \qquad\;\;\;
\end{array}
\right.
\ee
where 
\begin{eqnarray}\label{DS-1}
\begin{array}{c}
f_k : X_k\times \M_k\times \Y_k\rightarrow X_{k+1},\qquad\qquad\\
\\
X_k,\;X_N\subset \R^m,\;\;\;\;\M_k\subset \R^n,\qquad\qquad\\
\end{array}
\end{eqnarray}
for every $k=0,\ldots,N-1$. Here $x_k$ is the state space, $u_k$ the control and $y_k$ the random
disturbance. For every $k$, we are given the following constraint functions:
\begin{eqnarray}\label{c-f}
c_1^k,\ldots,c_{j_k}^k: {\cal M}_k\times A_k\rightarrow \R,\;\;j_k< n,
\end{eqnarray}
where $A_k$ is an open subset of $\R^m$, with $X_k\subset A_k$.\\
The set of admissible controls at time $t_k$ for the point $x\in X_k$ is defined as follows:
\begin{eqnarray}\label{c-U}
U_k(x)=\{u\in \M_k\;:\;c_i^k(u,x)\leq 0,\;i=1,\ldots,j_k\}.
\end{eqnarray}

We assume that $U_k(x)$ is non empty, for every $x\in X_k$ and that the random disturbance $y_k$ is
a measurable function over a Probability Space $(\Omega,\F,P)$ with values in a measurable space
$(\Y_k,{\cal E}_k)$, where ${\cal E}_k$ is a sigma-field over $\Y_k$.\\  
The disturbance $y_k$ is characterized by a probability law $p_k(\cdot)$, which we assume
independent of $(x_k,u_k)$ and of prior disturbances $y_{k-1},\ldots,y_0$.\\

We call an {\em admissible control law}, a set $\phi=\{u_0,\ldots,u_{N-1}\}$ of functions
$u_k:X_k\rightarrow \M_k$ such that $u_k(x)\in U_k(x)$, for every $x\in X_k$.\\
We denote by ${\cal U}$ the set of admissible control laws.\\

Given an initial state $x\in X_0$, the optimization problem consists in finding an admissible
control law $\phi \in {\cal U}$ which maximizes the cost functional
\begin{eqnarray}\label{OPT}
J_{\phi}(x)=\E[g(x_N)]
\end{eqnarray}
where $\E$ denotes the expected value taken over $(\Omega,\F, P)$, and $x_N$ is the value at time $t_N$
of the state $x$, according to (\ref{DS}).\\
The real-valued function $g:X_N\rightarrow \R$ is called $stopping-cost$ or $utility$ function.\\
We want to point out that, eventhough we limit ourselves to this simpler model for the sake of
simplicity, the discussion in this paper can be easily extended to include running costs or
probability measures $p_k$ depending also on the state and control variable.\\

The Classical Dynamic Programming Algorithm consists in solving the problem (\ref{OPT})
by means of the following sequence of one-step optimization problems:
\be
\label{DP}
\left\{
\begin{array}{lr}
J_{k}(x)=\invstackrel{\sup}{u\in U_k(x)}\E_k[J_{k+1}(f_k(x,u,y_k))], &
   \forall\ x\in X_k,\ \ 0\leq k<N\\ \\
J_N(x)= g(x),& \forall\ x\in X_N.
\end{array}
\right.
\ee
The value function at time $t_k$, $J_k$ is defined over the state space $X_k$,
and $\E_k$ denotes the expected value taken w.r.t. the measure $p_k$ over $\Y_k$.
For any given initial state $x_0\in X_0$, the value $J_0(x_0)$ computed by the algorithm equals
the optimal cost
\begin{eqnarray}\label{DP-1}
\invstackrel{\max}{\phi \in {\cal U}}\  J_{\phi}(x_0)
\end{eqnarray} 
and the optimal control policy $\phi^*$ can be obtained by $\phi^* = \{u_0^*,\ldots,u_{N-1}^*\}$ 
where $u_k^*$, $k=0,\ldots,N-1$ maximizes the right-hand side of (\ref{DP}). See \cite{Be}
for a detailed description of the Dynamic Programming Algorithm in the discrete-time case.

%%%%%%%%%%%%%%%%%%%%%%%%%%%%%%%%%%%%%%%%%%%%%%%%%%%%%%%%%%%%%%%%%%%%%%%%%%%%%%%%%%%%%%%%%%%
\section{The Hausdorff Metric}\label{HAUSDORFF}
%%%%%%%%%%%%%%%%%%%%%%%%%%%%%%%%%%%%%%%%%%%%%%%%%%%%%%%%%%%%%%%%%%%%%%%%%%%%%%%%%%%%%%%%%%%

We introduce in this section the Hausdorff metric over the space of all compact subsets of a given
metric space. This metric is used in Section \ref{LIPSCHITZEST} to estimate the Lipschitz constant
of the multifunction $x\mapsto U_k(x)$.\\

To study the regularity of the value function $J_k$ defined in (\ref{DP}), since the admissible
control set $U_k(x)$ depends on the state variable, one needs to measure the distance between
the admissible control sets corresponding to different states of the system. \\
We need, therefore, to introduce a distance between sets in order to show some regularity
property of the function
$$
A \mapsto \max_{A} f
$$
with $f$ continuous and $A$ a subset of a separable metric space $(M,d)$.
Since we only consider compact control sets, we can limit ourselves to introducing the Hausdorff
metric on the class of all compact subsets of $M$, denoted by $Comp(M)$.\\
For every $K_1,\;K_2\in Comp(M)\backslash\{\emptyset\}$, let
\begin{eqnarray}\label{34}
d_H(K_1,K_2)=\inf\big\{\var>0:\;K_1\subset K^{\var}_2,\; and\; K_2\subset K_1^{\var} \big\}
\end{eqnarray}
where, for any set $A\subset X$, $A^{\var}=\{y:\;d(y,A)<\var\}$ denotes the open ball of radius
$\var$ around $A$. \\
It is easy to verify that $d_H$ is a metric on $Comp(M)$. Furthermore if two points $x,y$ of $X$ are
regarded as the single point sets $\{x\}$ and $\{y\}$ in $Comp(M)$, then
\begin{eqnarray}\label{34.1}
d_H(\{x\},\{y\})=d(x,y).
\end{eqnarray}
That is to say, $M$ is isometrically embedded in $Comp(M)$. See \cite{BARNSLEY} for the properties of 
the Hausdorff distance.\\

The following Proposition concerns the Lipschitz regularity of a real valued map defined over the
space of all compact subsets of a metric space. 
 
\begin{prop}\label{LIPFTILDE}
Let $f$ be a Lipschitz continuous function over the space $M$, and
\begin{eqnarray}\label{34.2}
\begin{array}{c}
\t{f}:Comp(M)\rightarrow \R.\\
K\mapsto \invstackrel{\max}{K} f.
\end{array}
\end{eqnarray}
Then $\t{f}$ is a Lipschitz continuous map over the metric space $(Comp(M),d_H)$ and its Lipschitz constant equals 
the Lipschitz constant of $f$ over $M$.
\end{prop}

\begin{Dimo}{\bf Proposition \ref{LIPFTILDE}.}
Let $K_1,K_2$ be compact subsets of $M$, then there exist $x_i\in K_i$, such that
$\t{f}(K_i)=f(x_i)$, for $i=1,2$. By the definition $(\ref{34})$, for a fixed $\delta>0$, there exists $\var>0$, such that
\begin{eqnarray}\label{34.3}
\var&<&d_H(K_1,K_2)+\delta,\nonumber\\
&& K_1\subset K_2^{\var},\nonumber\\
&& K_2\subset K_1^{\var}.
\end{eqnarray} 
Therefore there exist $y_i\in\K_i$ $i=1,2$, such that,
\begin{eqnarray}\label{34.4}
d(x_1,y_2)&<&d_H(K_1,K_2)+\delta,\nonumber\\
d(y_1,x_2)&<&d_H(K_1,K_2)+\delta.\nonumber\\
\end{eqnarray}
Hence, by the definition of $x_1,x_2$ and the Lipschitz regularity of $f$ it follows,
\begin{eqnarray}\label{34.5}
\t{f}(K_1)-\t{f}(K_2)&\leq& f(x_1)-f(y_2)\leq Lip(f)d(x_1,y_2)\nonumber\\
&\leq& Lip(f)( d_H(K_1,K_2)+\delta),\nonumber\\
\t{f}(K_1)-\t{f}(K_2)&\geq& f(y_1)-f(x_2)\geq -Lip(f)d(y_1,x_2)\nonumber\\
&\geq& -Lip(f)( d_H(K_1,K_2)+\delta).\nonumber\\
\end{eqnarray}
The previous inequalities, being $\delta$ arbitrary, proves
\begin{eqnarray}\label{34.6}
|\t{f}(K_1)-\t{f}(K_2)|\leq Lip(f) d_H(K_1,K_2).
\end{eqnarray}
Therefore $\t{f}$ is a Lipschitz continuous map over $Comp(M)$, and $Lip(\t{f})\leq Lip(f)$.\\
On the contrary, if $x_1,x_2\in X$, then by (\ref{34.1}), we have 
\begin{eqnarray}\label{34.7}
|f(x_1)-f(x_2)|=|\t{f}(\{x_1\})-\t{f}(\{x_2\})|\leq Lip(\t{f})d(x_1,x_2)
\end{eqnarray}
that proves $Lip(f)\leq Lip(\t{f})$.
\cvd
\end{Dimo}

%%%%%%%%%%%%%%%%%%%%%%%%%%%%%%%%%%%%%%%%%%%%%%%%%%%%%%%%%%%%%%%%%%%%%%%%%%%%%%%%%%%%%%%%%%%
\section{Main Notations and Definitions}\label{notations}
%%%%%%%%%%%%%%%%%%%%%%%%%%%%%%%%%%%%%%%%%%%%%%%%%%%%%%%%%%%%%%%%%%%%%%%%%%%%%%%%%%%%%%%%%%%

In this section we introduce the main notations which are used in this paper.\\
\ben
\item Let $f$ be some real valued, Lipschitz continuous function over the domain $\D\subset
\R^m$. We refer to the Lipschitz constant of $f$ as to
\begin{eqnarray}\label{7.3}
Lip(f)=\sup_{\begin{array}{c} 
   (x,x^\p)\in \D\times \D\\
    x\neq x^\p   
   \end{array}} \frac{|f(x)-f(x^\p)|}{|x-x^\p|}.
\end{eqnarray}
\item We recall the definition of norm of an operator $S:\R^{p_1}\rightarrow \R^{p_2}$
\begin{eqnarray*}\|S\|:=\invstackrel{\max}{\begin{array}{c}
x\in\R^{p_1}\\		
|x|=1\end{array}}|Sx|.\end{eqnarray*}
It follows :
\bar\label{disug}
\|S\|\leq \sqrt{p_1 p_2}\invstackrel{\max}{\begin{array}{c} h=1,\ldots,p_2\\l=1,\ldots,p_1\end{array}}|S_{h,l}|=:\sqrt{p_1 p_2}\;\|S\|_\infty.
\ear
\item If $c\in \R^j$ we say $c\leq 0$ if and only if $c_i \leq 0$ for every $i=1,\ldots,j$.
\item Let $d\geq j$ be two integers, and $\Pi$ be the set of all multi-indexes
$\pi=(i_1,\ldots,i_j)$, with $1\leq i_1<i_2<\ldots<i_j\leq d$.
Then for every $\pi\in \Pi$ and $u\in\R^d$, let $u_\pi$ denote the projection of $u$ over the coordinates specified by $\pi$, i.e.
$u_{\pi} = (u_{i_1},\dots,u_{i_j}) \in \R^j$.\\
For every differentiable function $f:V\times A\rightarrow \R^j$, with $V\subset \R^d$ and $A\subset\R^m$ open sets, and for any $\pi\in \Pi$, let $\frac{\partial f}{\partial v_\pi}$ denote
the Jacobian matrix of $f$ w.r.t. the coordinates of $v\in V$ specified by $\pi$, i.e.
$$ 
\left( \frac{\de f}{\de v_\pi} \right)_{hl} = \frac{\de f_h}{\de v_{i_l}}, \qquad h,l=1,\ldots,j 
$$
\item Given $\pi$ in $\Pi$, $d\geq j$, we define the map
\begin{eqnarray}\label{NOT}
Z^{\pi}:\R^d\times \R^j\rightarrow \R^d
\end{eqnarray}
where, for every $u\in\R^d$, $v\in\R^j$,
\begin{eqnarray}\label{NOT-1}
(Z^{\pi}(u,v))_\pi=v,\;\;\;Z^{\pi}_i(u,v)=u_i,\;\;\;\forall\;i\notin\pi,\;i=1,\ldots,d.
\end{eqnarray}
In other words $Z^{\pi}(u,v)$ is a obtained from $u$, by substituting the vector $v$ to the
components of $u$ corresponding to $\pi$. Obviously if $j=d$, we have the only
$\pi=(1,2,..,d)$ and $Z^{\pi}(u,v) = v$.
\item For every $\pi=(i_1,\ldots,i_j)\in \Pi$, let denote by $T_\pi$, the
matrix
\bar\label{AUX.9-1}
(T_\pi)_{h,l}=\left\{
\begin{array}{lr}
1 & h = i_l \\
0 & \mbox{otherwise}
\end{array}\right.
\ear
for every $1\leq h\leq d$, and $1\leq l\leq j$.
\item We call {\em regular arc} a function $\gamma:{\cal I}\rightarrow \R^m$, where ${\cal I}$
is a compact interval, which is piecewise differentiable on ${\cal I}$, with $\gamma^\p$ bounded and 
$\gamma^\p(t)\neq 0$, for every $t\in{\cal I}$ where the derivative exists.\\
\een
\begin{Def}\label{CONNECT}
A non empty subset $X\subset\R^m$ satisfies the property
($CON$), if $X$ is connected and there exists a characteristic number $a(X)>0$, such that
for every pair of distinct points $x_1,\;x_2\in X$, there exists a regular
arc $\gamma :[w_1,w_2]\rightarrow X$
such that the following conditions hold true:
\begin{eqnarray}\label{CON}
\gamma(w_i)&=&x_i,\;\;i=1,2,\nonumber\\ \nonumber\\
l(\gamma)&\leq &a(X)|x_1-x_2|,
\end{eqnarray}
where $l(\gamma):=\int_{w_1}^{w_2} |\gamma^\p(t)|dt$, denotes the length of $\gamma$.
\end{Def}
Obviously, if $X$ is a convex subset of $\R^m$, then the property $(CON)$ holds with $a(X)=1$.
Also, it can be proved that any connected compact submanifold of $\R^m$ satisfies the property
$(CON)$.\\

\begin{Def}\label{LipMAN}
A set $\M\subset\R^n$ is called a $d$-dimensional {\bf Lipschitz manifold} if, for every $u\in \M$,
there exists an open neighborhood $W$ of $u$ in $\M$ and an homeomorphism $\psi :
W\rightarrow V$, where $V$ is an open subset of $\R^d$ and such that $\psi$ and
$\psi^{-1}$ are Lipschitz continuous. The couple $(W,\psi)$ is called a {\bf local chart} for
$u$.
\end{Def}

\begin{Def}
An {\bf atlas} on a compact Lipschitz manifold $\M\subset\R^n$ is a finite collection of charts
$\{(W_{\alpha},\psi_{\alpha})\}_{\alpha\in A}$ such that $\{W_{\alpha}\}_{\alpha\in A}$ is a
covering of $\M$.
\end{Def}

\begin{Def}
\label{NRM}
Let $\M$ be a Lipschitz manifold, we define ${\cal NR}(\M)$ to be the set of nonregular points of
$\M$, that is, for every point $u\in \M\backslash{\cal NR}(\M)$ there exists a local chart
$(W,\psi)$ for $u$ whose inverse is $C^1$ in a neighborhood of $\psi(u)$.
Such a chart will be called a {\bf regular chart for $u$}.\\
\end{Def}

Up to a suitable change of the atlas, in the reminder of the paper we make the following assumption:
\begin{assumption}
\label{chart-ass}
If $u\in \M\backslash {\cal NR}(\M)$, every local chart $(W,\psi)$ for $u$ is regular.
\end{assumption}

\begin{oss}
With Assumption \ref{chart-ass}, the Definition \ref{NRM} becomes:
$$
\M \backslash {\cal NR}(\M) = \left\{ u\in \M :\ \mbox{every local chart}\ (W,\psi)\ \mbox{for}\
u\ \mbox{is a regular chart for}\ u \right\}
$$
\end{oss}

\begin{oss}
\label{varie-manifold}
Let $\M\subset \R^n$ be a Lipschitz manifold with dimension $d$.\\
\ben
\item We denote by $T_u(\M)$ the tangent space to $\M$ at $u\in \M$. $T_u(\M)$
exists $H_\M^d$-almost everywhere on $\M$ by Rademacher's Theorem, where $H_\M^d$ is the
$d$-dimensional Hausdorff measure on $\M$. Since $\M\subset \R^n$, we can view $T_u(\M)$ as a linear
subspace of $\R^n$. 
In particular, if $u\notin{\cal NR}(\M)$ and $\phi: V\subset \R^d \avalori \phi(V)\subset \M$,
$\phi(v) = u$, is the inverse of a regular chart for $u$, the $n$-dimensional vectors
\be
\label{a.1}
\frac{\partial \phi}{\partial v_1}(v),\ldots,\frac{\partial \phi}{\partial v_d}(v),
\ee 
form a basis of $T_u(\M)$.\\
\item A function $f:\M\rightarrow \R^j$ is differentiable at $u\notin {\cal NR}(\M)$ ($resp.$ $C^1$
$in$ $a$ $neighborhood$ $of$ $u$) if for every chart $\psi$ defined in a neighborhood of $u$ with an
inverse differentiable at $\psi(u)$ ($resp.$ $C^1$ $in$ $a$ $neighborhood$ $of$ $\psi(u)$), the
function $f\circ\psi^{-1}$ is differentiable at $\psi(u)$ ($resp.$ $C^1$ $in$ $a$ $neighborhood$
$of$ $\psi(u)$) in the usual sense.\\
\item Let $f:\M\rightarrow \R^j$ be differentiable at $u\notin {\cal NR}(\M)$. The {\bf differential
of $f$ at $u$} is a linear operator
\bar\label{differential}
d_u f: T_{u}(\M)\rightarrow \R^j
\ear
defined as follows. Let $\phi$ be the inverse of a regular chart for $u$, as in 1.,
and let be given the basis (\ref{a.1}) of the tangent space $T_u(\M)$. Then, for any $w\in T_u(\M)$,
$$
w=\sum_{i=1}^d w_i\frac{\partial \phi}{\partial v_i}(v),\;\;\;w_i\in\R,\;i=1,\ldots,d,
$$
The differential of $f$ at $u$ is defined as
\be
\label{diff}
d_u f(w)=\sum_{i=1}^d w_i\frac{\partial f\circ\phi}{\partial v_i}(v)\in\R^j.
\ee
It can be proved that this definition is well-posed in that it does not depend on the choice of the
chart $\phi^{-1}$.\\
Given the immersion of $T_u(\M)$ in $\R^n$, we can view the differential of a map as a linear
operator over a $d$-dimensional linear subspace of $\R^n$ and therefore we can consider its
norm. In other words $\M$ is a Riemannian manifold.\\
The same definition of differential holds if $f$ depends also on a variable $x\in A\subset\R^m$. In
this case definition (\ref{diff}) applies to $f(\cdot,x)$ for every $x\in A$, and its differential
is denoted by $d_u f(\cdot,x)$.\\
\een
\end{oss}

\begin{Def}
Let $\M \subset \R^n$ be a Lipschitz manifold and $f:\M\times A \avalori \R^j$, $A$ an open subset
of $\R^m$, $X\subset A$. We define ${\cal NR}(f)$ to be the set of nonregular points of $f$, i.e.
\be
\label{NRf} 
{\cal NR}(f) = \left\{ u\in \M\backslash {\cal NR}(\M):\ \exists\ x\in A:\ f\ \mbox{is
non-differentiable at}\ (u,x)\right\}.
\ee
\be
\label{NSf}
{\cal NS}(f,X) = \left\{ u\in \M\backslash \left({\cal NR}(\M) \unione {\cal NR}(f)\right):\
\exists\ x\in X:\ f(u,x) \leq 0\ \mbox{and}\ d_uf(\cdot,x)\ \mbox{is not surjective}\right\}
\ee
and
\be
\label{D}
\D(\M,X,f) = {\cal NR}(\M) \unione {\cal NR}(f) \unione {\cal NS}(f,X)
\ee
\end{Def}

\begin{Def}
\label{TuDef}
For every $u \in \M\backslash {\cal NR}(\M)$, we choose, once and for all, a regular chart for $u$,
whose inverse will be denoted by $\phi_u$, $\phi_u(v) = u$.
We define, for any $\pi \in \Pi$, the linear operator
\be
\label{Tu}
T^u_{\pi}:\R^j \avalori T_u(\M)
\ee
whose associated matrix, w.r.t. the canonical basis of $\R^j$ and the basis
$$
\frac{\partial \phi_u}{\partial v_1}(v),\ldots,\frac{\partial \phi_u}{\partial v_d}(v)
$$
of $T_u(\M)$, is the matrix $T_{\pi}$ defined in (\ref{AUX.9-1}).
\end{Def}

%%%%%%%%%%%%%%%%%%%%%%%%%%%%%%%%%%%%%%%%%%%%%%%%%%%%%%%%%%%%%%%%%%%%%%%%%%%%%%%%%%%%%%%%%%%
\section{The Lipschitz Regularity of the Value Function}\label{LIPSCHITZEST}
%%%%%%%%%%%%%%%%%%%%%%%%%%%%%%%%%%%%%%%%%%%%%%%%%%%%%%%%%%%%%%%%%%%%%%%%%%%%%%%%%%%%%%%%%%%

In order to prove that the value function is Lipschitz continuous we need to prove that the map
$U_k(\cdot)$ defined in (\ref{c-U}) is Lipschitz continuous with respect to the Hausdorff distance
$d_H$ in $Comp(X_k)$.\\
We present some regularity results for multifunctions which have the same structure as the set of
admissible controls (\ref{c-U}).\\

\begin{teo}\label{TECH-REG}
Let ${\bf c} = (c_1,\ldots,c_j)$, with $c_i$, $i=1,\ldots,j$, $j\leq n$, real valued maps defined
on $\M\times A$, where $\M\subset \R^n$ and $A\subset \R^m$ is open.  Let $X$ be a non empty subset
of $A$ which has the property ($CON$), and let $\M$ be a compact Lipschitz $d$-dimensional manifold,
with $j \leq d$. We further assume that
\begin{description}
\item[$i$)] For every $x\in X$ the set
\begin{eqnarray}\label{i}
U(x)= \{ u\in \M\;:\;{\bf c}(u,x)\leq 0 \}
\end{eqnarray}
is non empty.\\
\item[$ii$)] $\D(\M,X,{\bf c})$ is closed and one of the following assumptions holds:
\begin{description}
\item[{\bf (A)}] The set $\{(u,x)\;:\;x\in X,\;u\in U(x)\cap \D(\M,X,{\bf c})\}=\emptyset$.\\
\item[{\bf (B)}] $\D(\M,X,{\bf c})\neq \emptyset$ and every 
$u\in U(x)\cap \D(\M,X,{\bf c})$, with $x\in X$, is of adherence for $U(x)\backslash\D(\M,X,{\bf c})$.
\end{description}
\item[$iii$)] The function ${\bf c}:=(c_1,\ldots,c_j)$ is $C(\M\times A)\cap
C^1(\M\backslash\D(\M,X,{\bf c})\times A)$.\\
\item[$iv$)] For every $x \in X$ and $u \in U(x)\backslash\D(\M,X,{\bf c})$,
\be
\label{ii-2}
\tau:=\invstackrel{\sup}{x\in X}\; \invstackrel{\sup}{ u\in\; U(x)\backslash\D(\M,X,{\bf c})}\;{\cal T}(u,x)<\infty\;.
\ee
with
\be
\label{G(.)}
{\cal T}(u,x):=\invstackrel{\max}{\pi\in \Pi(u,x)}\left\|\left(d_{u} {\bf c}(\cdot,x)\circ
T^u_{\pi}\right)^{-1}\circ d_x c(u,\cdot)\right\|,
\ee
and
\be
\Pi(u,x):=\left\{\pi\in \Pi\;:\;d_{u} {\bf c}(\cdot,x)\circ T^u_{\pi}\;\mbox{ is invertible} \right\}.
\ee
\end{description}
Then the map
\begin{eqnarray}\label{TECH2}
x\in X\mapsto U(x)
\end{eqnarray}
is
\begin{description}
\item[if $ii$)-{\bf (A)} holds]: $d_H$-Lipschitz continuous and its Lipschitz constant can be
estimated by $a(X)\tau Lip_{\M} $, where $a(X)$ denotes the characteristic number of $X$ as in
Definition \ref{CONNECT} and
\bar\label{lipc}
Lip_{\M}:=\sup\left\{Lip(\phi):\;\phi^{-1}\mbox{ is a chart over $\M$}\right\}.
\ear
\item[if $ii$)-{\bf (B)} holds]: $d_H$-uniformly continuous on every compact subset of $X$ which has the
property ($CON$).
\end{description}
\end{teo}

Next result strenghtens the regularity assumptions on the constraint function ${\bf c}$ in order to obtain
Lipschitz regularity of the value function even in the case $ii$)-{\bf (B)}.\\

\begin{teo}\label{TECH-REG-2}
In the hypotheses $i)$,$ii)$-{\bf (B)}, $iii)$ of Theorem \ref{TECH-REG}, we assume that
\begin{description} 
\item[$iv$)] there exist $\mu, \lambda, r > 0$ such that for every $x\in X$ and $u \in
U(x)\backslash\D(\M,X,{\bf c})$, and $\pi \in \Pi(u,x)$,
\bar
\label{M(.)}
\left\{\begin{array}{l}
\left\|\left(d_{u} {\bf c}(\cdot,x)\circ T^u_{\pi}\right)^{-1}\right\|\leq \l,\\
\\
\|d_y c(u,\cdot)\|\leq \mu,\ \forall\ y\in A, |y-x| \leq r
\end{array}\right.
\ear
\end{description}
and that
\begin{description} 
\item[$v$)] if $r$ is chosen as in $iv$), there exists $L>0$ such that for every $u\in
U(x)\backslash\D(\M,X,{\bf c})$ with $x\in X$, and for every regular chart $(W,\psi)$ for $u$ it
holds:
\be
\label{second}
\left\| \frac{\partial \varphi}{\partial v}(v,y) - \frac{\partial \varphi}{\partial v}
(v_u,x) \right\| \leq L \left( \left| v - v_u \right| + \left| y - x \right| \right)
\ee
for every $v$ in a suitable neighborhood of $v_u$, $y\in A$, $|y-x| \leq r$, with
$\varphi(\cdot,\cdot) = {\bf c}(\psi^{-1}(\cdot),\cdot)$, $v_u = \psi(u)$.
\end{description}
Let 
\be
\label{tau2}
\tau = \lambda\mu,
\ee
then $x\in X\mapsto U(x)$ is $d_H$-Lipschitz continuous, with constant $a(X)\tau Lip_\M$.
\end{teo}

\begin{oss}
\label{Remark4.7}
If {\bf (A)} holds true, then for every compact set $K\subset X$, we have
$$
\{(u,x)\;:\;x\in K,\;u\in U(x)\cap \D(\M,X,{\bf c})\}=\emptyset,
$$
hence, by compactness, there exists $\s^*>0$ such that 
$U(x)\cap \left(\D(\M,X,{\bf c})\right)^{\s^*}=\emptyset$, for any $x\in K$.
\end{oss}

Using these results, we prove a regularity result for the value function (\ref{DP}). Then we
prove Theorem \ref{TECH-REG} and Theorem \ref{TECH-REG-2}.\\

\begin{teo}\label{LIP}
Let be given (\ref{DS})-(\ref{c-U}) and the related optimization algorithm 
(\ref{DP}). Suppose that:
\begin{description}
\item[1)] for every $k=0,\ldots,N-1$, the triplet $(\M_k,{\bf c}^k, X_k)$ satisfies the assumptions
of Theorem \ref{TECH-REG} under the condition $ii)$-{\bf (A)}, or the hypotheses of Theorem
\ref{TECH-REG-2}.\\
\item[2)] for every $k=0,\ldots,N-1$, there exists a nonnegative, ${\cal E}_k$-measurable and $p_k$-integrable function $y\in \Y_k\mapsto V_k(y)$, i.e.:
\begin{eqnarray}\label{p-law}
\int_{\Y_k}V_k(y)dp_k(y)<\infty,
\end{eqnarray} 
such that for $p_k-a.e.\; y\in\Y_k$,
\begin{eqnarray}\label{f_k}
\begin{array}{c}
|f_k(x,u,y)-f_k(x^\p,u^\p,y)|\leq V_k(y)|(x-x^\p,u-u^\p)|,\qquad \qquad \qquad\\
\\
\;\;\forall\;\;x,x^{\p}\in X_k,\;\;\;\forall\;\; u,u^{\p}\in \M_k.\qquad \qquad \qquad
\end{array}
\end{eqnarray}
\end{description}
If $g$ is Lipschitz continuous over $X_N$ then, for every $k$, $J_k$ is Lipschitz continuous over $X_k$
and the following estimate holds:
\be
\label{LIP0}
\left\{
\begin{array}{c}
Lip(J_k) \leq L(J_{k+1}) \E_k[V_k](1+a_k\tau_k Lip_{\M_k}),\;\;\;k=0,\ldots,N-1\\
\\
Lip(J_N) =Lip(g).\;\;\;\;\;\;\;\;\;\;\; \qquad\qquad\qquad\qquad\qquad\qquad\qquad\qquad\\
\end{array}
\right.
\ee
with $a_k$ the characteristic number of $X_k$, given in Definition \ref{CONNECT}, 
and $\tau_k$ defined in (\ref{ii-2}) or (\ref{tau2}) according to the two alternatives of assumption
{\bf 1)}.
\end{teo}

We prove this result first, then we prove Theorems \ref{TECH-REG} and \ref{TECH-REG-2}.\\

\begin{Dimo}{\bf Theorem \ref{LIP}.}  
We proceed by induction over $k$. For $k=N$, the function $J_N$ is the utility function $g$, which
is supposed Lipschitz continuous over $X_N$. We assume $J_{k+1}$ Lipschitz continuous over $X_{k+1}$,
for $k\leq N-1$. Therefore, the function
$$
\Psi_k(x,u)=\E\left[J_{k+1}(f_k(x,u,y_k))\right],\;\;\;\;(x,u)\in X_k\times \M_k.
$$
is Lipschitz continuous over $X_k\times \M_k$, and an easy computation yields
\begin{eqnarray}\label{LIP3}
Lip(\Psi_k)\leq Lip(J_{k+1})\E_k[V_k].
\end{eqnarray} 
Moreover, for any $k=0,\ldots,N-1$, the assumption {\bf 1)} implies that the set-valued map $x\in
X_k\mapsto U_k(x)$ is $d_H$-Lipschitz continuous, by Theorem \ref{TECH-REG} or Theorem
\ref{TECH-REG-2}, with Lipschitz constant estimated by $a_k \tau_k Lip_{\M_k}$. Hence we can apply
Proposition \ref{LIPFTILDE}, to deduce that the marginal function defined on $X_k$ by
$$
x\in X_k\mapsto\invstackrel{\max}{u\in U_k(x)}\Psi_k(x,u)
$$
is Lipschitz continuous. The DP-algorithm (\ref{DP}), implies that this map is exactly $J_k$.
Moreover, by Proposition 2.1,  its Lipschitz constant is estimated by $Lip(\Psi_k) (1+Lip(U_k))$. 
Therefore we have
$$
Lip(J_k)\leq Lip(J_{k+1})\E_k[V_k]\left(1+a_k\tau_k Lip_{\M_k}\right),
$$
which proves the assertion.
\cvd
\end{Dimo}

In order to prove Theorem \ref{TECH-REG}, we have to estimate the $d_H$-distance between $U(x)$ and
$U(y)$ i.e. we need to prove, by the definition of $d_H$, an inclusion of the type $U(x) \subset (U(y))^{\delta}$, where $\delta>0$.\\
We need therefore to take an admissible control $u$ for $x$ and to show that we can ``perturb'' it to
an admissible control for $y$ whose distance from $u$ is less than $\delta$. 
The idea of the proof is to show that this property holds true locally and that the radius of the 
neighborhood where the property holds true is independent of the point $x$ and of the control $u$. This
independence will allow us to prove Theorem \ref{TECH-REG} by ``iteration'', i.e. by covering the arc 
between $x$ and $y$ with a finite number of balls of constant radius where the property holds true.

\begin{lem}\label{TECH-2}
In the same hypotheses of Theorem \ref{TECH-REG} and according to the two alternatives of
assumption $ii$), we have:
\begin{description}
\item[{\bf (A)}] for every regular arc $\gamma:[w_1,w_2]\rightarrow X$, there exists $r_0>0$ such that, for every $w_1 \leq t < s\leq w_2$ which satisfy
\begin{eqnarray}\label{TTT}
\gamma([t,s])\subset B_{r_0}(\gamma(t)),
\end{eqnarray}
the following inclusion holds:
\begin{eqnarray}\label{TECH-PROOF-7}
U(\gamma(t))\subset (U(\gamma(s)))^{\tau^{\p}},\;\;\;\forall\;\tau^{\p}>\tau_{t,s} Lip_\M,
\end{eqnarray}
with
\bar\label{length}
\tau_{t,s}:=\tau l(\gamma; t,s):=\tau \int_t^s |\gamma^\p(\xi)|d\xi, 
\ear
and $\tau$ is given by (\ref{ii-2}).
\item[{\bf (B)}] For every compact set $K\subset X$ which has the property ($CON$), and for any
$\s>0$, there exists $r_\s\in (0,\s]$ such that, for every $x,y\in K$ such that
\bar\label{TTT-U}
|x-y|<r_\s,
\ear
we have
\begin{eqnarray}\label{TECH-PROOF-7-U}
U(x)\subset (U(y))^{ \tau^\p+\s},\;\;\;\forall\;\tau^{\p}> a(K)\tau Lip_\M|x-y|.  
\end{eqnarray}
with $a(K)$ the characteristic constant for $K$, given in Definition \ref{CONNECT}.
\end{description}
\end{lem}

Let's assume for the time being that Lemma \ref{TECH-2} holds true. Using this result it is
straightforward proving Theorem \ref{TECH-REG}:\\

\begin{Dimo} {\bf Theorem \ref{TECH-REG}.}
If $X$ is a single point set, then we do not need to prove anything. Let's assume that $X$ contains
at least a pair of distinct points $x_1$, $x_2$; we want to estimate the Hasudorff distance between
the corresponding admissible control sets, i.e.
\begin{eqnarray}\label{m5}
d_H\left(U(x_1),U(x_2)\right).
\end{eqnarray}
By the property ($CON$), there exists a regular arc $\gamma:[w_1,w_2]\rightarrow X$ which connects
$x_1$ to $x_2$ and which satisfies (\ref{CON}). We have to distinguish between the two alternatives
in assumption $ii$)-{\bf(B)}:\\

{\bf Case (A):} For every $t\in [w_1,w_2]$ let denote by $B_t$ the open ball centered in $\gamma(t)$
with radius $r_0>0$, obtained by applying Lemma \ref{TECH-2}, case {\bf (A)}, to the arc
$\gamma$. We can assume that the curve does not have self-intersections, otherwise we could define a
new arc connecting the same points and having smaller length. We introduce the following sequence
\begin{eqnarray}\label{m1}
\left\{\begin{array}{l}
t_0:=w_1\;\;\;\;\;\;\;\;\;\;\;\;\;\;\;\;\;\;\;\;\;\;\;\;\;\;\;\;\;\;\;\;\;\;\;\;\;\;\;\;\;\;\;\;\;\;\;\;\;\;\;\;\;\;\;\;\;\;\;\;\;\;\;\;\;\;\;\;\;\;\;\;\;\;\;\;\;\;\;\;\;\;\;\;\;\;\;\;\\
\\
t_{i+1}:=\sup\{w_2\geq t\geq t_i:\;|\gamma(t_i)-\gamma(s)|\leq \frac{r_0}{2},\;\forall\;s\in [t_i,t]\},\;\;\mbox{ if $i\geq 0$,}
\end{array}\right.
\end{eqnarray}
that is increasing in the interval $[w_1,w_2]$. It is straightforward proving that $t_p=w_2$ for
some integer $p>0$.\\
Therefore $\{ [t_0,t_1],\ldots,[t_{p-1},t_p]\}$ is a finite covering of $[w_1,w_2]$. Moreover
$\gamma([t_i,t_{i+1}])\subset B_{t_i}$.\\
By (\ref{TECH-PROOF-7}), we obtain
\begin{eqnarray}\label{m3}
U(\gamma(t_i))\subset \left(U(\gamma(t_{i+1}))\right)^{\tau_i},\;\;\;i=0,\ldots,p-1,
\end{eqnarray}
where $\tau_i=\tau Lip_\M\int_{t_i}^{t_{i+1}}|\gamma(\xi)|d\xi+\frac{\var}{p}$, $\var$ is an arbitrary positive number. 
By iterating the previous inclusion (\ref{m3}), we obtain
\begin{eqnarray}\label{m4}
U(x_1)\subset \left(U(x_2)\right)^{\tau Lip_\M l(\gamma)+ \var}.
\end{eqnarray}
By switching the role of $x_1$ and $x_2$ and by (\ref{CON}), the definition of the metric $d_H$ and
the arbitrary choice of $\var$, we can estimate (\ref{m5}) with $a(X)\tau Lip_\M|x_1-x_2|$.\\

{\bf Case (B):} Let $K$ be a compact subset of $X$ which satisfies the condition ($CON$).
By applying Lemma \ref{TECH-2} in the case {\bf (B)}, for any $\s>0$, there exists $r_\s\in
(0,\s]$, such that the inclusion (\ref{TECH-PROOF-7-U}) holds, whenever $x,y\in K$ satisfy
(\ref{TTT-U}). This implies
$$
d_H(U(x),U(y))\leq (1+a(K)\tau Lip_\M)\s,\;\;\;\forall\;\;|x-y|<r_\s,
$$
and therefore $U(\cdot)$ is uniformly continuous over $K$, with respect to $d_H$.
\cvd
\end{Dimo} 

We turn now to the proof of Lemma \ref{TECH-2}. In the case ({\bf B}), we need to approximate the
controls which lie in a neighborhood of $\D(\M,X,{\bf c})$.

\begin{lem}\label{AUX}
Let $\M$ be a manifold as in Theorem \ref{TECH-REG}, $K\subset X$ compact, and suppose that
assumptions $i$), $ii$)-{\bf (B)} and $iii$) hold true.\\ 
For every $\s>0$, there exists $\rho>0$ such that for every $x\in K$ and $u\in U(x)$, there exists
$u_\s\in U(x)$, such that
\bar\label{AUX.1}
d(u_\s,\D(\M,X,{\bf c}))\geq \rho,\;\;\;\mbox{   and   }\;\;\;|u_\s-u|\leq \s.
\ear
\end{lem}

\begin{Dim}{\bf .}
We prove the result for the controls $u\in U(x)\intersez(\D(\M,X,{\bf c}))^{\rho}$, otherwise it suffices
choosing $u_{\sigma} = u$.\\
Suppose by contradiction that there exists $\o{\s}>0$ such that for $\rho=\frac{1}{h}$, there exist
$u_h\in U(x_h)$, $x_h\in K$ with
\bar\label{AUX.3}
d(u_h,\D(\M,X,{\bf c}))<\frac{1}{h}
\ear
such that, for any $u\in U(x_h)$,
\bar\label{AUX.4}
d(u,\D(\M,X,{\bf c}))<\frac{1}{h},\;\;\;\mbox{ or }\;\;\;|u-u_h|>\o{\s}.
\ear
Without loss of generality we may assume that $x_h\rightarrow \o{x}$, $u_h\rightarrow \o{u}$, as
$h\rightarrow \infty$, for some $\o{x}\in K$ and $\o{u}\in \D(\M,X,{\bf c})$, since $K$ is compact and
$\D(\M,X,{\bf c})$ is closed by assumption $ii$) of Theorem \ref{TECH-REG}. \\
The continuity of ${\bf c}$ over $\M\times A$, implies
$\o{u}\in U(\o{x})\cap \D(\M,X,{\bf c})$ and, by $ii$), $\o{u}$ is of adherence for
$U(\o{x})\backslash\D(\M,X,{\bf c})$, therefore there exists $\t{u}\in U(\o{x})\backslash\D(\M,X,{\bf c})$ such that
\bar\label{AUX.5}
|\t{u}-\o{u}|<\frac{\o{\s}}{2}.
\ear
Let $\pi \in \Pi(\tilde{u},x)$ and $\phi = \phi_{\tilde{u}}$ be the map, relative to $\tilde{u}$,
fixed in Definition \ref{TuDef}, with $\phi(\t{v})=\t{u}$, $\t{v}\in\R^d$.
Let $Z^\pi$ be the function introduced in Section \ref{notations}, and
\bar\label{AUX.6}
F(w,y):={\bf c}\left(\phi(Z^\pi(\t{v},w)),y\right)-{\bf c}(\t{u},\o{x})\in\R^j,\;\;\;\forall\;(w,y)\in \R^j\times A,
\ear
which is of class $C^1$ in a neighborhood of $((\tilde{v})_\pi, \ov{x})$. We observe that
\bar\label{AUX.7}
F((\t{v})_{\pi},\o{x}) = {\bf c}\left(\phi(\t{v}), \o{x}\right) - c(\t{u},\o{x}) = 0
\ear
and
\bar\label{AUX.8}
\frac{\de F}{\de w}\left((\t{v})_{\pi},\o{x}\right) = \frac{\de {\bf c}(\phi(\cdot),\o{x})}{\de v}(\t{v},\o{x})\cdot \frac{\de Z^{\pi}}{\de w}, 
\ear
where the first matrix on the right-hand side is $j\times d$ and the second one is $d\times j$. It
is easy to verify that if $\pi = (i_1,\ldots,i_j)$, $1\leq i_1<i_2<\ldots<i_j\leq d$, and if $1\leq
i \leq d$, $1\leq l \leq j$, we have
\bar\label{AUX.9}
(T_\pi)_{i,l}=\left(\frac{\de Z^{\pi}}{\de w}\right)_{i,l} = \left\{
\begin{array}{lr}
1 & i = i_l \\
0 & \mbox{otherwise}
\end{array}
\right.
\ear
By the definition of $\pi$ we deduce the invertibility of the Jacobian matrix (\ref{AUX.8}).
With (\ref{AUX.7}) this allows for the application of the classical Implicit Function Theorem which
implies the existence of a continuous map $q:B\rightarrow \R^j$, where $B\subset A$ is an open
neighborhood of $\o{x}$ in $\R^m$, such that
\bar\label{AUX.10}
q(\o{x})=(\t{v})_\pi
\ear
and
\bar\label{AUX.11}
F(q(y),y)=0,\;\;\;\forall\;y\in B.
\ear
Since ${\bf c}(\t{u},\o{x})\leq 0$, using (\ref{AUX.6}) and (\ref{AUX.11}), we get
\bar\label{AUX.12}
\t{q}(y):=\phi(Z^{\pi}(\t{v},q(y)))\in U(y),\;\;\forall\;y\in B.
\ear
For large values of $h$, we have $x_h\in B$ and since $\t{q}(x_h)\rightarrow \t{u}$ and
$u_h\rightarrow \o{u}$ as $h\rightarrow \infty$, by (\ref{AUX.5}), we also may assume  
$$
|\t{q}(x_h)-u_h|<\o{\s}.
$$
Hence by (\ref{AUX.12}) and (\ref{AUX.4}) with $u=\t{q}(x_h)$, we infer
\bar\label{AUX.13}
d(\t{q}(x_h),\D(\M,X,{\bf c}))<\frac{1}{h},
\ear
for large $h$. Letting $h\rightarrow \infty$ in (\ref{AUX.13}),
we obtain
$$
\t{u}=\t{q}(\o{x})\in \D(\M,X,{\bf c}),
$$
which is a contradiction.
\cvd
\end{Dim}

The main difficulty in the proof of Lemma \ref{TECH-2} is in building the radius $r_0$ independent
of $t,s$ and of the controls in $U(\gamma(t))$. We use a quantitative version of the classical
Implicit Function Theorem, which provides an estimate of the radius of the balls where the implicit
map is defined, see \cite{Chierchia}.\\
Using the implicit function theorem, we are able to build a map that, for every state
$\gamma(s)$ ``near'' $\gamma(t)$, prescribes how to modify the control $u\in U(\gamma(t))$ to obtain
an admissible control for $\gamma(s)$.

\begin{teo}\label{IFT}
Let $F: {\cal O}\rightarrow \R^j $ be a map defined in the open 
set ${\cal O}\subset\R^j\times \R^m$. Let $r_1$, $r_2>0$
be such that, if $B_1=\{v\in\R^j:\;|v-v_0|\leq r_1\}$ and $B_2=\{y\in\R^m:\;|y-y_0|\leq r_2\}$, 
then
$B_1\times B_2\subset {\cal O}$ and the following hypotheses hold:
\bar\label{reg}
F,\;\;\frac{\partial F}{\partial v}\;\;\mbox{ {\em are continuous on $B_1\times B_2$}},
\ear
\begin{eqnarray}\label{IFT-1}
F(v_0,y_0)=0,\qquad\qquad\qquad\det\frac{\partial F}{\partial v}(v_0,y_0)\neq 0,
\end{eqnarray}
\begin{eqnarray}\label{IFT-4}
\left\{\begin{array}{c}
\invstackrel{\sup}{y\in B_2}|F(v_0,y)|\leq\frac{r_1}{2\|T_0\|},\qquad \qquad T_0=\left(\frac{\partial F}{\partial v}(v_0,y_0)\right)^{-1}\\
\\
\invstackrel{\sup}{B_1\times B_2}\left\|I_j-T_0\frac{\partial F}{\partial v}\right\|\leq \frac{1}{2},\qquad\qquad\qquad\;\;\;\qquad\qquad
\end{array}\right.
\end{eqnarray}
where $I_j$ is the identity matrix of order $j$. Then there exists a unique
function $q\in C(B_2;B_1)$ which satisfies
\begin{eqnarray}\label{IFT-2}
q(y_0) =v_0,
\end{eqnarray}
and for every $(v,y)\in B_1\times B_2$, it holds
\begin{eqnarray}\label{IFT-3}
F(v,y)=0\;\;\; \Longleftrightarrow \;\;\;v=q(y).
\end{eqnarray}
\end{teo}
Next Proposition provides an expression for the Jacobian matrix of the implicit function $q$.

\begin{prop}\label{IFT-reg}
Let $D_1$, $D_2$ be two open balls of $\R^j$ and $\R^m$, respectively.
Suppose that $F\in C^1(D_1\times D_2; \R^j)$, and $q\in C(D_2;D_1)$ satisfies $F(q(y),y)=0$, for every $y\in D_2$.
If $\frac{\partial F}{\partial v}$ is invertible in $D_1\times D_2$, then $q\in C^1(D_2)$ and
\begin{eqnarray}\label{
IFT-5}
\frac{\partial q}{\partial y}(y)=-\left(\frac{\partial F}{\partial v}(q(y),y)\right)^{-1}\frac{\partial F}{\partial y}(q(y),y),\mbox{ for $y\in D_2$ }.
\end{eqnarray}
\end{prop}

We turn finally to the proof of Lemma \ref{TECH-2}.\\

\begin{Dimo}{\bf Lemma \ref{TECH-2}.} We prove the assertions in steps.
We derive an approximation of $U(x)$ by the controls of $U(y)$ for $x,y$ which lie in 
a compact connected subset $K$ of $X$, then we specialize the discussion according to
the assumptions $ii$)-${\bf (A)}$ and $ii$)-{\bf (B)}.

{\bf (Construction of a covering).} Let $K$ be a nonempty compact connected subset of $X$, $\o{\s}>0$ such that
\bar\label{variab}
K^{\o{\s}}\subset A.
\ear
Let $\s^*$ be chosen as in Remark 5.8 and $\s<\o{\s}$. We introduce
\bar\label{sigma}
\rho:=\left\{\begin{array}{c}
\s^*\;\;\;\mbox{ if {\bf (A)} holds}\\
\\
\rho_\s\;\;\;\mbox{ if {\bf (B)} holds,}\\
\end{array}\right.
\ear
where $\rho_\s$ is related to $\s$ and to the compact set $K\subset X$ via Lemma \ref{AUX}.
Let 
\begin{eqnarray}\label{m7}
K_{\gamma,\rho}=\{(u,x)\;:\;x\in K,\;\;u\in U(x)\mbox{  and  }d(u,\D(\M,X,{\bf c}))\geq \rho\}.
\end{eqnarray} 
By Remark \ref{Remark4.7}, in the case {\bf (A)}, and Lemma \ref{AUX}, in the case {\bf (B)}, for
small $\s$, this set is a non empty, compact subset of $\M\times K$.
For every $(u,x)\in K_{\gamma,\rho}$, let $\pi(u,x) \in \Pi(u,x)$ and $\phi_u$ the map given
in Definition \ref{TuDef}, for the point $u$, $\phi_u:V_u \avalori \M$. We may assume, without loss of
generality, that $\phi \in C^1(V_u)$, furthermore, since $\D(\M,X,{\bf c})$ is closed, we may also
suppose $\phi_u(V_u) \intersez \D(\M,X,{\bf c}) = \emptyset$.\\
By the assumption $iii$), the map $\varphi_u$, defined by 
$$\varphi_u(v,y)={\bf c}(\phi_u(v),y),\;\;\;\;\forall\;(v,y)\in V_u\times A,$$
is $C^1$. 
Furthermore, the matrix
\be
\label{Malpha}
R(w,y;u,x) := \frac{\de \varphi_{u}}{\de v_{\pi(u,x)}}
\left(\phi_u^{-1}(w),y\right),\;\;\;\forall\;(w,y)\in\phi_u(V_u)\times A,
\ee
is invertible at $(u,x)$. \\
By the continuity of the function (\ref{Malpha}), there exists $\delta^\p=\delta^\p(u,x)>0$ such that
\bar\label{m8}
z(u,x):=\inf\left\{\left|\det\;R(w,y;u,x)\right|\;:\;(w,y)\in (B_{\delta^\p}(u)\cap\M)\times B_{\delta^\p}(x) \right\}>0,
\ear
where $B_{\delta^\p}(u)\subset\R^n$ is the open ball of radius $\delta^\p$ centered at $u$, with
$B_{\delta^\p}(u)\cap\M\subset \phi_u(V_u)$, and $B_{\delta^\p}(x)\subset K^{\o{\s}/2}$ is the open ball in $\R^m$ of radius $\delta^\p$ centered at $x$. 
By the continuity of $\phi_u^{-1}$ there exists $\delta(u,x)\leq \delta^\p(u,x)$ such that
\bar
\label{neigh}
|\phi^{-1}_u(w)-\phi^{-1}_u(u)|\leq \frac{\zeta_u}{2},\;\;\;\forall\;\;w\in B_{\delta(u,x)}(u)\cap \M.
\ear
with $\zeta_u > 0$ such that the closed ball in $\R^d$ centered at $\phi_u^{-1}(u)$ with radius
$\zeta_u$ is contained in $V_u$. Such a ball is denoted by ${\cal N}_u$.\\
The collection
$$ \left\{B_{\frac{\delta}{2}}(u)\times B_{\frac{\delta}{2}}(x)\; :\; \delta=\delta(u,x),\;\;(u,x)\in K_{\gamma,\rho} \right\} $$
is an open covering of $K_{\gamma,\rho}$, therefore we can extract a finite covering corresponding to some
points $(u_1,x_1),\ldots,(u_p,x_p)$ $\in K_{\gamma,\rho}$.
To simplify the remainder of the proof let us define:
$$
\begin{array}{l}
\phi_i:=\phi_{u_i},\;\;\varphi_i:=\varphi_{u_i},\;\;R_i(\cdot,\cdot):=R(\cdot,\cdot;u_i,x_i),\;\;T_i:=T^{u_i}
_{\pi_i},\;\;L_i:=Lip(\phi_i),\\
\\
\pi_i:=\pi(u_i,x_i),\;\;\;\;z_i := z(u_i,x_i),\;\;\;\;\delta_i := \delta(u_i,x_i),\;\;\;\; \zeta_i :=
\zeta_{u_i}\\ \\
\N_i:=\N_{u_i},\;\;\;\,O_{i}(\delta) := (B_{\delta}(u_i) \intersez \M)\times B_{\delta}(x_i), \qquad
\forall \delta > 0,\;\;\;\;i=1,\ldots,p, \\
\end{array}
$$
and
\be
\begin{array}{l}
\delta^*:=\invstackrel{\min}{1\leq i\leq p} \delta_i,\;\;\;z^* := \invstackrel{\min}{1\leq i\leq p}
z_i,\;\;\;\;\theta^{*}:=\invstackrel{\min}{1\leq i\leq
p}\frac{\delta_i}{L_i},\;\;\;\;\;\;\;\;\zeta^* := \invstackrel{\min}{1\leq i \leq p}\ \zeta_i
\end{array}
\ee

By (\ref{m8}), we have
\bar\label{stronzo1}
\left| \det R_{i}(w,y) \right| \geq z^* > 0,\;\;\;\;\forall\;(w,y)\in O_{i}(\delta_i),\;\;\forall\;i=1,\ldots,p. 
\ear
Now suppose $(u,x)\in K_{\gamma,\rho}$, then $(u,x)\in O_i(\frac{\delta_i}{2})$, for some $i=1,\ldots,p$. Let 
$$
\begin{array}{c}
v_i:= \phi_i^{-1}(u) \in \R^d,
\end{array}
$$
and
$$
\begin{array}{c}
F:\R^j\times A \avalori \R^j
\end{array}
$$
given by
\be
\label{F}
F(v,y) := \varphi_i\left(Z^{\pi_i}(v_i,v),y\right) - c(u,x)
\ee
\\
which is defined and $C^1$ in the open ball in $\R^d$ centered at $(v_i)_{\pi_i}$ and with radius
$\frac{\zeta_i}{2}$, since (\ref{neigh}) implies 
$$
\left| Z^{\pi_i}(v_i,v) - \phi_i^{-1}(u_i) \right| < \zeta_i
$$
that implies $Z^{\pi_i}(v_i, v) \in V_{u_i}$.\\
In order to apply Theorem \ref{IFT} to $F$ with $(v_0,y_0) =
((v_i)_{\pi_i}, x)$
, let us assume for the time being that $r_1, r_2$ are chosen as prescribed by
Theorem \ref{IFT}. We observe that $F(v_0,y_0) = 0$ and 
$$
\frac{\de F}{\de v}\left((v_i)_{\pi_i}, x\right) = R_{i}(u,x),
$$
is invertible by (\ref{stronzo1}). By Theorem \ref{IFT}, there exists a continuous map
$$ q:B_2 \avalori B_1 $$
with $B_1\times B_2 \subset \R^j\times A$, such that (\ref{IFT-2}) and (\ref{IFT-3}) hold true. The
function $q$ is a feedback function that allows us to build an admissible control for $y$ starting
from the admissible control $u$ for $x$ : as in the proof of Lemma \ref{AUX} (see (\ref{AUX.12})),
we have
\bar\label{ffff-REG}
\t{q}(y):=\phi_{i}\left(Z^{\pi_i}(v_{i}, q(y))\right) \in U(y)\backslash\D(\M,X,{\bf c}),
\ear
for every $y\in B_2$.\\

{\bf (Construction of $r_1$ and $r_2$).} In order to construct $r_1, r_2$ as in Theorem \ref{IFT}, we consider
\bar\label{lambda}
\l :=\invstackrel{\max}{1\leq l\leq p}\;\;\invstackrel{\sup}{(w,y)\in O_{l}(\frac{\delta_l}{2})}\left\|\left(R_l(w,y)\right)^{-1}\right\|,
\ear
which is finite since for every $l$
$$ \left\| \left( R_{l}(w,y) \right)^{-1}\right\| \leq \frac{Const}{z^*},$$
and $Const$ is a constant that depends only on the supremum of the norm of $\frac{\partial
\varphi_l}{\partial v}$ over a compact subset of $\R^d\times\R^m$.
The regularity assumption $iii$) allows us to define
\begin{eqnarray}\label{TECH-PROOF-4}
\mu &>& \sup\left\{\left\|\frac{\partial {\bf c }}{\partial x}(u^\p,x^\p)\right\|\;:\;u^\p\in
\M,\;\;d(u^\p,\D(\M,X,{\bf c}))\geq \rho,\;\;x^\p\in K^{\o{\s}/2}\right\}\\
\label{TECH-aaaa}
\omega(h,k)&:=&\invstackrel{\max}{1\leq l\leq p}\ \sup\Big\{\left\|\frac{\de \varphi_l}{\de v}(v,y)-\frac{\de \varphi_l}{\de v}(v^\p,y^\p)\right\|\;:\;|v-v^\p|\leq h,\;|y-y^\p|\leq k,\nonumber\\
&& v,\;v^\p\in {\cal N}_{l},\;y,\;y^\p\in K^{\o{\s}/2}\Big\},\;\;\;\;\;h,\;k\geq 0.
\end{eqnarray}
Let us fix
\be
\label{condr1}
\left\{\begin{array}{l}
r_2 = \beta r_1,\;\;\;\beta = \frac{1}{2 \mu \lambda},\\
\\
r_1 < \min\left(\frac{\zeta^*}{2},\frac{\o{\s}}{2\beta},\frac{\delta^*}{2\beta},\frac{\theta^*}{4}\right).
\end{array}\right.
\ee
Notice that the modulus $\omega$ defined in (\ref{TECH-aaaa}) is a descreasing function of $h$ and
$k$ and its limit for $h,k\rightarrow 0$ is zero, so we may choose $r_1$ so that
\bar\label{cccccc}
\omega(r_1,\beta r_1)\leq \frac{1}{2\l\sqrt{dj}}.
\ear
Let us verify the inequalities in (\ref{IFT-4}) with $(v_0, y_0) = ((v_{i})_{\pi_i}, x)$:
\bar\label{1-dis}
\invstackrel{\sup}{y\in B_2}|F(v_0,y)|=\invstackrel{\sup}{y\in B_2}|{\bf c}(u,y)-{\bf c}(u,x)|\leq \mu r_2=\frac{r_1}{2\l}\leq \frac{r_1}{2\|T_0\|},\;\;T_0=\left(R_i(u,x)\right)^{-1}
\ear
where we used $d(u,\D(\M,X,{\bf c}))\geq \rho$ and $x+\eta(y-x)\in K^{\o{\s}/2}$, for every $y\in B_2$ and
$0\leq \eta\leq 1$. Furthermore for any $(v,y)\in B_1\times B_2$ we have
\bar\label{2-dis}
\!\!\left\|I-T_0\frac{\partial F}{\partial v}(v,y)\right\|\!\!\! &\leq& \!\!\!\|T_0\|\left\|\frac{\partial  \varphi_i}{\partial v}(v_i,x)-\frac{\partial \varphi_i}{\partial v}(Z^{\pi_i}(v_i,v),y)\right\|\|T_i\|\nonumber\\
&\leq & \l\sqrt{dj} \omega(r_1,\beta r_1)\leq \frac{1}{2}.
\ear
This inequality follows by the definition of $\omega$ and by (\ref{cccccc}): in fact
$v_i,\;Z^{\pi_i}(v_i,v)\in {\cal N}_i$, and $r_2<\frac{\o{\s}}{2}$ implies $x,y\in K^{\o{\s}/2}$. 
This proves (\ref{IFT-4}) and justifies the 
application of Theorem \ref{IFT}.\\
The radius $r_2$ depends only on $K$ and $\s$, and it does not depend on the
particular choice of $(u,x)$.\\

{\bf (Approximation of the control $u$).} Again we need to distinguish between the alternatives {\bf
(A)} and {\bf (B)} in $ii$):\\

{\bf Case (A).} Let $K$ be the image of the arc $\gamma$ and $r_0=r_2(K,\s^*)$.
We consider $s>t$ as in (\ref{TTT}), where $x=\gamma(t)$, and we
have to approximate $u$ by an admissible control for the state $\gamma(s)$. To this purpose we
apply Proposition \ref{IFT-reg} to the pair $F$, $q$ obtained in the previous step. 
With the choice (\ref{condr1}) for $r_1$, let
\bar\label{D_1D_2}
\begin{array}{c}
D_1:=\left\{v\in\R^j\;:\;|v-(v_{i})_{\pi_i}|<2r_1\right\}\\
\\
D_2\;\mbox{ is the interior of $B_2$}.\;\;\;\;\;\;\;\;\;\;\;\;\;\;\;\;\;\;
\end{array}
\ear
For any $v\in D_1$ and $y\in D_2$,
$$
\begin{array}{c}
|\phi_i(Z^{\pi_i}(v_i,v))-u_i|\leq |\phi_i(Z^{\pi_i}(v_i,v))-\phi_i(v_i)|+|u-u_i|< 2 L_i r_1 +\frac{\delta_i}{2}\leq \delta_i,\\
\\
|y-x_i|\leq  r_2+|x-x_i|<\frac{\delta^*}{2}+\frac{\delta_i}{2}\leq \delta_i.\;\;\;\;\;\;\;\;\;\;\;\;\;\;\;\;\;\;\;\;\;\;\;\;\;\;\;\;\;\;\;\;\;\;\;\;\;\;\;\;\;\;\;\;\;\;\;\;\;\;\;\;\;\;\;\;\;
\end{array}
$$
This implies 
\bar\label{ddddd1}
(\phi_i(Z^{\pi_i}(v_i,v)),y)\in O_i(\delta_i),\;\;\;\forall\;(v,y)\in D_1\times D_2.
\ear
Therefore by (\ref{stronzo1}), the matrix
\bar\label{p-F}
\frac{\partial F}{\partial v}(v,y)=R_i(\phi_i(Z^{\pi_i}(v_i,v)),y)
\ear
is invertible over $D_1\times D_2$, and $q(D_2)\subset B_1\subset D_1$. 
By applying Proposition \ref{IFT-reg} we deduce that $q\in C^1(D_2)$. Since 
$\gamma([t,s])\subset D_2$, for any $\xi\in [s,t]$, by (\ref{ffff-REG}) 
\bar\label{me}
\t{q}(\gamma(\xi))\in U(\gamma(\xi))\backslash\D(\M,X,{\bf c}).
\ear
The assumption $iv$) implies
\bar\label{fff}
\left\|\frac{\partial q}{\partial y}(\gamma(\xi))\right\|&\leq&
\left\| [R_i(\t{q}(\gamma(\xi)),\gamma(\xi))]^{-1}\frac{\partial {\bf c}}{\partial x}(\t{q}(\gamma(\xi)),\gamma(\xi))\right\|\nonumber\\
&=&\left\|\left(d_{\t{q}(\gamma(\xi))}{\bf c}(\cdot,\gamma(\xi))\circ T^{\t{q}(\gamma(\xi))}_{\pi_i}\right)^{-1}\circ d_{\gamma(\xi)}{\bf c}(\t{q}(\gamma(\xi),\cdot)\right\|\nonumber\\
\nonumber\\
&\leq& {\cal T}(\t{q}(\gamma(\xi)),\gamma(\xi))\leq \tau.
\ear
In fact the matrix associated to the linear operator 
$d_{\t{q}(\gamma(\xi))}{\bf c}(\cdot,\gamma(\xi))\circ T^{\t{q}(\gamma(\xi))}_{\pi_i}$, with respect
to the canonical basis of $\R^j$, is $R_i(\t{q}(\gamma(\xi)),\gamma(\xi))$.	
Using (\ref{fff}) we obtain
\bar\label{TECH-PROOF-6}
|\t{q}(\gamma(s))-u| &\leq&  Lip_{\M}|q(\gamma(s)) - q(\gamma(t))|\nonumber\\                      
&\leq & \tau Lip_{\M}\int_t^s|\gamma^\p(\xi)|d\xi \leq \tau_{t,s}Lip_{\M}.
\ear
This inequality proves that
\bar\label{same}
u\in \left(U(\gamma(s))\right)^{\tau^\p},\;\;\;\forall\;\tau^\p>Lip_{\M}\tau_{t,s}.
\ear
Since $r_2$ is independent of $(u,x)$, we have
\begin{eqnarray}\label{TECH-PROOF-7-1}
U(\gamma(t))\backslash\left(\D(\M,X,{\bf c}))\right)^{\rho}\subset\left(U(\gamma(s))\right)^{\tau^\p}.
\ear

{\bf Case (B).} Let $K$ be a compact subset of $X$, which has the property ($CON$), and
$r_\s=\min(\s,\frac{r_2(K,\s)}{a(K)})$. Let $x,y\in K$ be such that $|x-y|<r_\s$. Then we consider a
regular arc $\gamma:[w_1,w_2]\rightarrow K$ as in the Definition \ref{CONNECT}, which connects $x$
to $y$ and lies in $K$. For every $\xi\in [w_1,w_2]$, we have
$$
|\gamma(\xi)-\gamma(w_1)|\leq l(\gamma)\leq a(K)|x_1-x_2|<r_2,
$$
which implies
$$
\gamma([w_1,w_2])\subset B_{r_2}(\gamma(w_1)).
$$
Hence by the same arguments developed in the previous case, we deduce
\begin{eqnarray}\label{TECH-PROOF-7-2}
U(x)\backslash\left(\D(\M,X,{\bf c}))\right)^{\rho}\subset\left(U(y)\right)^{\tau^\p},\;\;\;\forall\;\tau^\p>a(K)\tau Lip_\M|x_1-x_2|.
\ear

{\bf (Conclusions).} 
\begin{description}
\item[{\bf (A)}] Using the inclusion (\ref{TECH-PROOF-7-1})
and the definition of $\rho$ in (\ref{sigma}), we infer
\begin{eqnarray}\label{TECH-PROOF-7-A}
U(\gamma(t))\subset\left(U(\gamma(s))\right)^{\tau^\p}.
\ear
In fact, in this case, $\rho=\s^*$, which is defined in Remark 5.8 so that
$U(x)\subset (\D(\M,X,{\bf c}))^{\s^*}$
for every $x\in X$.
\item[{\bf (B)}] Let $x,y\in K$, $u\in U(x)$, with $d(u,\D(\M,X,{\bf c}))<\rho=\rho_\s$. By Lemma \ref{AUX}, we find
$u_\s\in U(x)$, with $d(u_\s,\D(\M,X,{\bf c}))\geq \rho$, such that
$$
|u-u_\s|\leq \s.
$$
Hence, using (\ref{TECH-PROOF-7-2}), there exists $v_\s\in U(y)$ such that $|u-v_\s|<\tau^\p+\s$, $\tau^\p>a(K)Lip_\M \tau$.
This implies
\bar\label{TECH-PROOF-7-B}
U(x)\subset\left(U(y)\right)^{\tau^\p+\s}.
\ear
\end{description}
This concludes the proof of Lemma \ref{TECH-2}.
\cvd
\end{Dimo}

\begin{Dimo}{\bf Theorem \ref{TECH-REG-2}.}
Let, for any $u\in \M$ and for any local chart $(W,\psi)$ for $u$,
\be
\label{deltauwpsi}
\delta(u,W,\psi) = \sup \left\{ \delta > 0:\ \o{B_{\delta}(\psi(u))} \subset \psi(W) \right\}
\ee
and
\be
\label{deltau}
\delta(u) = \invstackrel{\max}{(W,\psi):\ u\in W}\ \delta(u,W,\psi).
\ee
Since every local chart $\psi$ is Lipschitz continuous and $\M$ is compact in $\R^n$, $\delta(u)$ is
finite for every $u\in \M$. We prove that the map $\delta(\cdot)$ is uniformly bounded from
below. Since $\M$ is compact, it suffices proving that $\delta(\cdot)$ is lower semi-continuous,
i.e. that the set 
\be
\label{set1}
\left\{ u\in \M:\ \delta(u) > \alpha \right\}
\ee
is open in $\M$, for any $\alpha >
0$. Suppose, for the time being, that the set
\be
\label{set2}
\left\{ u\in W:\ \delta(u,W,\psi) > \alpha \right\}
\ee
is open in $\M$ for any local chart $(W,\psi)$ in the atlas of $\M$ and for any $\alpha > 0$. Then
if $u$ is such that $\delta(u) > \alpha$, there exists a local chart $(W,\psi)$ for $u$ such that
$\delta(u,W,\psi) > \alpha$. Therefore, there exists $\eta > 0$ such that $\delta(w,W,\psi) >
\alpha$ for any $w\in B_{\eta}(u) \intersez W$. Hence $\delta(w) > \alpha$ for any $w \in
B_{\eta}(u) \intersez W$.  Since $W$ is open in $\M$, $B_{\eta}(u) \intersez W$ is open in $\M$, and
therefore $\delta(\cdot)$ is lower semi-continous over $\M$.\\
To prove that the set (\ref{set2}) is open in $\M$, let $u\in W$ be such that $\delta(u,W,\psi) > \alpha$
and $\xi \in (0, \frac{\delta(u,W,\psi) - \alpha}{Lip(\psi)})$ and $\eta \in (\alpha, \delta(u,W,\psi) -
Lip(\psi)\xi)$, then
$$
\o{B_{\eta}(\psi(w))} \subset \psi(W) \qquad \forall\ w\in B_{\xi}(u)\intersez W.
$$
This implies
$$
\delta(w,W,\psi) \geq \eta > \alpha \qquad \forall\ w\in B_{\xi}(u)\intersez W.
$$
Therefore, since $W$ is open in $\M$, the assertion is proved.\\

Let 
\be
\label{udelta}
\underline{\delta} = \invstackrel{\min}{u\in \M}\ \delta(u) > 0
\ee
and $\Gamma$ be the image of a regular arc $\gamma :[w_1,w_2]\rightarrow X$, which 
connects two distinct points $x_1,x_2\in X$ and such that the property ($CON$) holds. 
Let $\o{\s}>0$ be such that
$$
\Gamma^{\o{\s}}\subset A.
$$
Let $x\in \Gamma$ and $u\in U(x)\backslash\D(\M,X,{\bf c})$, $\pi \in \Pi(u,x)$ and
$\phi_u: V\rightarrow \M$, the map fixed for $u$ in
Definition \ref{TuDef}. We may assume, up to modifying from the beginning the choice of the local
charts $\{ \phi_u \}_{u\in M}$ that appear in Definition \ref{TuDef}, that $\phi_u$ is such that
$$
B_{\underline{\delta}}(\phi_u^{-1}(u)) \subset V.
$$
Let $\phi = \phi_u$, $W = \phi(V)$ and $F$ the map
\bar\label{AUX.6-REG}
F(w,y):=\varphi\left(Z^\pi(v_u,w),y\right)-{\bf c}(u,x)\in\R^j,\;\;\;\forall\;(w,y): Z^\pi(v_u,w)\in V, y\in A,
\ear
with $\varphi(\cdot,\cdot):={\bf c}(\phi(\cdot),\cdot)\in C^1(V\times A)$ and $\phi(v_u)=u$, then
\bar\label{IFT-rg}
\frac{\partial F}{\partial w}(w,y) = \frac{\partial \varphi}{\partial v_\pi}(Z^\pi(v_u,w),y).
\ear
Since $F((v_u)_\pi,x)=0$ and using (\ref{IFT-rg}), the assumptions (\ref{IFT-1}) are
satisfied. Therefore we can apply Theorem \ref{IFT}, with $(v_0,y_0):=((v_u)_\pi,x)$. 
The inequalities in (\ref{IFT-4}) can be proved as follows, by assumptions $iv)$ and $v)$: let $r_1$
and $r_2$ be given by
\be
\left\{
\begin{array}{l}
r_2 = \beta r_1, \qquad \beta = \frac{1}{2\lambda\mu}\\ \\
r_1 < \min \left(\frac{r}{2\beta}, \frac{\ov{\s}}{4\beta}, \frac{1}{2}\underline{\delta}, \frac{1}{4L 
(1+\beta)\sqrt{dj}\lambda}\right)\\
\end{array}
\right.
\ee
then since the matrix associated to the linear operator $d_u {\bf c}(\cdot,x)\circ T^u_\pi$,
w.r.t. the canonical base of $\R^j$, is
\be
T_0^{-1} = \frac{\partial \varphi}{\partial v_{\pi}}(v_u, x),
\ee
we obtain
$$ 
\invstackrel{\sup}{y\in B_2} |F(v_0,y)| = \invstackrel{\sup}{y\in B_2} |c(u,y) - c(u,x)| \leq $$
$$ 
\invstackrel{\sup}{y\in B_2} \| d_y c(u,\cdot) \| r_2 \leq \mu r_2 = \frac{r_1}{2 \lambda} \leq
\frac{r_1}{2 \|T_0\|}.
$$
As to the second inequality in (\ref{IFT-4}), since
$$
\left| Z^\pi(v_u, w) - v_u \right| \leq r_1 \qquad \forall\ w\in B_1
$$
and
$$
|y-x| \leq r_2 < r \qquad \forall\ y\in B_2
$$
we have
$$
\|I_j-T_0\frac{\partial F}{\partial w}(w,y)\| \leq \lambda \sqrt(dj) \left\| \frac{\de \varphi}{\de
v}(v_u, x) - \frac{\de \varphi}{\de v}(Z^\pi(v_u,w),y) \right| \leq \lambda \sqrt{dj} L (r_1 +
r_2) \leq \frac{1}{2}.
$$
Moreover, by Proposition \ref{IFT-reg}, the implicit function $q:D_2\rightarrow D_1$ is $C^1(D_2)$,
where $D_1$ and $D_2$ are defined as in (\ref{D_1D_2}). Let $t>s$, such that $x=\gamma(t)$, and
$\gamma([t,s])\subset D_2$, we derive the inequality (\ref{TECH-PROOF-6}), as in the proof of Lemma
\ref{TECH-2}.
Since $r_2$ and $r_1$ do not depend on $u$, we deduce the inclusion
$$
U(\gamma(t))\backslash\D(\M,\Gamma,{\bf c})\subset\left(U(\gamma(s))\right)^{\tau^\p},\;\;\;\forall\;\tau^\p>\tau_{t,s}Lip_\M.
$$
Taking the closure, by assumption {\bf (B)}, we conclude that
$$
U(\gamma(t))\subset\left(U(\gamma(s))\right)^{\tau^\p},\;\;\;\forall\;\tau^\p>\tau_{t,s}Lip_\M.
$$
In the case {\bf (A)}, we can repeat the argument used in the proof of Theorem \ref{TECH-REG}
substituting $r_2$ for $r_0$. This brings to the inequality (\ref{m4}) which, as before, implies
$$
d_H(U(x_1),U(x_2))\leq a(X)\tau Lip_\M|x_1-x_2|.
$$
This proves the Lipschitz regularity of $U(\cdot)$.
\cvd 
\end{Dimo}

%%%%%%%%%%%%%%%%%%%%%%%%%%%%%%%%%%%%%%%%%%%%%%%%%%%%%%%%%%%%%%%%%%%%%%%%%%%%%%%%%%%%%%%%%%%
\section{An Application to a Financial Problem}\label{application}
%%%%%%%%%%%%%%%%%%%%%%%%%%%%%%%%%%%%%%%%%%%%%%%%%%%%%%%%%%%%%%%%%%%%%%%%%%%%%%%%%%%%%%%%%%%

We discuss here an application of the results in the previous sections to an optimal
asset-allocation problem.
More precisely we consider an optimal asset-liability management model in presence of
constraints:  the company can manage the investment coming from the policy-holders' payments, in
order to satisfy several regulatory and solvency constraints and to achieve a given objective.\\
The company can decide, at each time step, how to distribute the total wealth between the available
assets to achieve its goal, but he has to obey a number of constraints.  Furthermore, each portfolio
adjustment entails transaction costs, since it means selling part of an asset to provide either
liquidity or a different asset (see \cite{nostro}, \cite{ALTROLAVORO} for a detailed description of the
financial model).\\

We assume, for the sake of simplicity, that the manager can choose at each time step
$t_0,\ldots,t_{N-1}$, between a riskless and a risky investment, denoted by $B$ and $S$,
respectively, though the procedure described applies also to the more general case of $n$ possible
investments characterized by different values of yield and volatility.\\
We assume that the transaction costs of moving wealth between the sections are paid only on buying
and not on selling and that these transaction costs are linearly proportional to the size of the
transaction.
The evolution equations for the amounts invested in stocks and bonds are
\begin{eqnarray}\label{model.1}
\left\{
\begin{array}{c}
S_{k+1}=[S_k(1-u_k)+B_k v_k(1-\l^s)]Y_k^s,\;\;\;S_0=s_0>0\\
\\
B_{k+1}=[B_k(1-v_k)+S_k u_k(1-\l^b)]Y_k^b,\;\;\;B_0=b_0>0,
\end{array}
\right.
\end{eqnarray}
where for every $0\leq k\leq N-1$, $(u_k,v_k)\in \M_k \subset [0,1]\times [0,1]$ represent the
percentage of risky investment that is moved to riskless investment and viceversa, at time $t_k =
k\Delta t$, $Y_k^s,Y_k^b :\Omega\rightarrow (0,\infty)$ are random variables and $0<\l^b,\l^s<1$
represent the transaction cost coefficients.\\

The model is self-financed in that it does not require the use of cash to perform the
transactions. After each portfolio adjustment, the investments evolve according to the stochastic
yields $Y_K^s, Y_k^b$. We assume that the joint process $\{(Y_k^s,Y_k^b)\}_k$ is a discrete Markov
chain which takes values in a sequence of finite discrete subspaces of $(0,\infty)^2$,
$\{\Y_k\}_k$. The chain is characterized by the density function:
\begin{eqnarray}\label{model.2}
p_k(y^s,y^b) = P(Y_k^s = y^s,Y_k^b = y^b),\;\;\;\forall\;(y^s,y^b)\in\Y_k.
\end{eqnarray}

The company has to satisfy at each time-step a regulatory constraint that imposes a limit on the
percentage of wealth invested in risky assets. More precisely, the adjustment must be such that the
fund after the adjustment satisfies:
\be
\frac{B_k(1-v_k) + S_k u_k(1-\lambda^b)}{B_k(1-v_k) + S_k u_k(1-\lambda^b) + S_k(1-u_k) + B_k
v_k(1-\lambda^s)} \geq \alpha
\ee
i.e. the percentage of wealth invested in riskless assets, after the portfolio adjustment, is bigger
than $\alpha$. Furthermore, since $\Y_k$ is finite, we may require that the adjustment is such that
the constraint is satisfied also at time $t_{k+1}$,
\be
\frac{B_{k+1}}{B_{k+1} + S_{k+1}} \geq \alpha
\ee
with $(S_{k+1},B_{k+1})$ given by (\ref{model.1}).

It can be shown that this constraint is equivalent to the following inequality:
\bar\label{model.5}
c_k(u_k,v_k,S_k,B_k):=S_k(1-u_k)+v_kB_k(1-\l^s)-q_k[B_k(1-v_k)+u_kS_k(1-\l^b)]\leq 0.
\ear
with 
$$
q_k = \frac{1-\alpha}{\alpha} \invstackrel{\min}{} \left(1, \invstackrel{\min}{(y^s,y^b)\in\Y_k}\frac{y^b}{y^s}\right).
$$

We would like to avoid that any of the two investments becomes null at some time. To this aim, we
observe that, because of the structure (\ref{model.1}), if the initial allocation $(s_0,b_0)$ is
such that $(s_0, b_0) \in X_0$, with
\be
\label{X0}
X_0=\{(S,B)\;:\;S,B\geq \Delta_0,\;S+B\leq D_0\},
\ee 
with $D_0>\Delta_0>0$ it is possible to define $\Delta_k$, $D_k$ and $\delta_k$ such that if the
control space at time $t_k$ is
\be
\label{Mk}
\M_k = [0,\delta_k]\times[0,\delta_k],
\ee
then $(S_k,B_k) \in X_k$ with
\bar
\label{xk}
X_k:=\left\{(S,B)\;:\;S,B\geq \Delta_k,\;S+B\leq D_k\right\}
\ear
and the admissible control space at time $t_k$, 
\be
\label{model.6}
U_k(S,B):=\{(u,v)\in \M_k : c_k(u,v,S,B)\leq 0\}
\ee
is nonempty for any $(S,B) \in X_k$.

By the structure (\ref{model.5}) of the constraints, if $U_k(S,B) \neq \emptyset$, then the control
$(\delta, 0)$ must belong to $U_k(S,B)$ for some $\delta \in (0,1)$. Therefore, let $\delta_k \in (0,1)$
be such that
$$
c_k(\delta_k,0,S,B) = S(1-\delta_k)-q_k[B+\delta_k S(1-\l^b)]\leq 0,\;\;\;\forall\;S,\;B>0,
$$
and therefore
\bar\label{deltak}
\delta_k=\invstackrel{\sup}{S,B>0}\frac{S-q_kB}{S(q_k(1-\l^b)+1)}=\frac{1}{q_k(1-\l^b)+1}\in (0,1),\;\;\;\forall\;k\geq 0.
\ear

Assuming $(S_k,B_k) \in X_k$ we can recursively compute $\Delta_{k+1}$ and $D_{k+1}$ that define
$X_{k+1}$. Since $S_k,B_k\geq \Delta_k>0$ and $u,v\in\M_k$, we have
$$
\begin{array}{c}
S_{k+1}\geq \Delta_k(1-\delta_k)\u{y}_k\\ \\
B_{k+1}\geq \Delta_k(1-\delta_k)\u{y}_k.
\end{array}
$$
where 
$$
\begin{array}{c}
\u{y}_k:=\min\left\{\min(y^s,y^b)\;:\;\;(y^s,y^b)\in \Y_k\right\}
\end{array}
$$
Therefore we define
$$
\Delta_{k+1}:=\Delta_k(1-\delta_k)\min(\o{y}_k^s,\o{y}_k^b).
$$
To obtain $D_{k+1}$ as a function of $D_k$, we observe that by equations (\ref{model.1}), if
$S_k+B_k\leq D_k$, for every possible choice of the controls in $\M_k$, we have
$$
S_{k+1}+B_{k+1}\leq D_k\o{y}_k,
$$
where 
$$
\begin{array}{c}
\o{y}_k:=\max\left\{\max(y^s,y^b)\;:\;\;(y^s,y^b)\in \Y_k\right\},
\end{array}
$$
that is
$$
D_{k+1}=D_k\o{y}_k,\;\;\;\forall\;k\geq 0.
$$

We have proved that if $(s_0, b_0) \in X_0$, with $X_0$ as in (\ref{X0}), equations (\ref{model.1})
map $X_k\times \M_k\times\Y_k$ into $X_{k+1}$, being $X_k$ defined as in (\ref{xk}). Furthermore the
admissible control set defined in (\ref{model.6}) is non empty for every $(S,B) \in X_k$.

As in section \ref{MODEL}, we want to find an optimal investment strategy which maximizes the
expected value of a given utility function $g$ at time $t_N$.\\

Since $\M_k=[0,\delta_k]\times[0,\delta_k]$ is not a Lipschitz manifold in the sense specified by
Definition \ref{LipMAN}, we decompose it as the union of infinitely many subsets, which
are $1$-dimensional Lipschitz manifolds. We can apply Theorem \ref{TECH-REG-2} to these manifolds
and then extend the result to $\M_k$.

\begin{prop}\label{TECH.model}
The map (\ref{model.6}) is $d_H$-Lipschitz continuous over $X_k$, for every $k\geq 0$. Moreover, if
$g:[0,+\infty)\times [0,+\infty) \avalori \R$ is Lipschitz continuous, the value function
$J_k(S,B)$, obtained via the DP algorithm (\ref{DP}) is Lipschitz continuous over $X_k$, for any
$k=0,\ldots,N$.
\end{prop}

\begin{Dim}{\bf.} 
We consider, for every $0\leq\delta< \delta_k$, the rectangle $R_{\delta,k}=[\delta,\delta_k]\times [0,\delta_k]$, and we define
\bar\label{delta-man}
\M_{\delta,k}=\partial R_{\delta,k},
\ear
and
\bar\label{delta-U}
U_{\delta,k}(S,B):=U_{k}(S,B)\cap \M_{\delta,k}
\ear
for every $(S,B)\in X_k$. We prove, by applying Theorem \ref{TECH-REG-2}, that the multifunctions
(\ref{delta-U}) are Lipschitz continuous over $X_k$ and that their Lipschitz constant does not depend
on the parameter $\delta$.\\
Since the space $X_k$ is convex, $a_k=a(X_k)=1$; furthermore $\M_{\delta,k}$ is obviously a
Lipschitz manifold in $\R^2$ of dimension $d=1=j$.\\

We want to apply Theorem \ref{TECH-REG-2} to derive the Lipschitz regularity of $U_{\delta,k}$.
The assumption $i$) in Theorem \ref{TECH-REG} is a consequence of the fact that $U_{\delta,k}(S,B)$ 
contains the control $(\delta_k,0)$, for every $(S,B)\in X_k$ and $\delta$.
Now let us observe that
\be
\label{ddddd}
\D(\M_{\delta,k},X_k,c_k)=\{(\delta,0),(\delta_k,0),(\delta_k,\delta_k),(\delta,\delta_k)\}.
\ee
Let $(S,B)\in X_k$ and $(u,v)\in \M_{\delta,k}\backslash{\cal NR}(\M_{\delta,k})$, the inverse of a
local chart for $(u,v)$ is 
\bar\label{delta-chart}
\phi(t):=\left\{\begin{array}{cc}
(t,0) & \mbox{ $\delta<t<\delta_k$, if $v=0$}\\
(\delta_k,t) & \mbox{ $0<t<\delta_k$, if $u=\delta_k$}\\
(t,\delta_k) & \mbox{ $\delta<t<\delta_k$, if $v=\delta_k$}\\
(\delta,t) & \mbox{ $0<t<\delta_k$, if $u=\delta$}.
\end{array}\right.
\ear
The constraint function $c_k$ is regular (${\cal NR}(c_k) = \emptyset$) and the jacobian matrix at
$(u,v)$ of the restriction of $c_k$ to the manifold $\M_{\delta,k}$ is,
by (\ref{model.5}),
\bar\label{delta-jacob}
J c_k(u,v)=\left\{\begin{array}{cc}
-S[1+q_k(1-\l^b)] & \mbox{ if $v=0$ or $v=\delta_k$}\\
B(1-\l^s+q_k) & \mbox{ if $u=\delta$ or $u=\delta_k$}.
\end{array}\right.
\ear
Therefore $\Pi(u,v,S,B)$ is non empty. This proves the assertion (\ref{ddddd}).

Now, $\D(\M_{\delta,k},X_k,c_k)$ is closed and, if
$(u,v)\in\D(\M_{\delta,k},X_k,c_k)\backslash\{(\delta_k,0)\}$ is an admissible control for $(S,B)$, it can
be approximated by
$$
(u_\var,v_\var)=\left\{\begin{array}{cc}
(\delta+\var,0) & \mbox{ if $(u,v)=(\delta,0)$}\\
(\delta_k,\delta_k-\var) & \mbox{ if $(u,v)=(\delta_k,\delta_k)$}\\
(\delta+\var,\delta_k) & \mbox{ if $(u,v)=(\delta,\delta_k)$}
\end{array}\right.
$$
for $0<\var<\delta_k,\delta_k-\delta$. In fact 
$(u_\var,v_\var)\in U_{\delta,k}(S,B)\backslash\D(\M_{\delta,k},X_k,c_k)$, 
as follows by the monotonicity of the constraint function $c_k$, w.r.t. $u$ and $v$.
If $u=\delta_k$ and $v=0$, then for every $(S,B)\in X_k$, by (\ref{deltak}) it holds
$$
c_k(u,v,S,B)=-q_k B\leq -q_k\Delta_k<0,
$$
therefore, by continuity, we can approximate $(\delta_k,0)$ with admissible controls of 
$U_{\delta,k}(S,B)\backslash\D(\M_{\delta,k},X_k,c_k)$. In other words the assumption $ii)$-{\bf
(B)} of Theorem \ref{TECH-REG} holds true.\\

To prove that the assumptions $iv),v)$ of Theorem \ref{TECH-REG-2} hold true, we observe that
\bar\label{delta-der-s-b}
\frac{\partial c_k}{\partial S}\in [0,1],\;\;\;\frac{\partial c_k}{\partial B}\in [-q_k,-q_k+\delta_k(q_k+1-\l^s)].
\ear
Therefore, the assumption $iv)$ follows by choosing
\bar\label{delta-M(.)}
\l_k:=\min\left(\frac{\delta_k}{\Delta_k},\frac{1}{\Delta_k (1-\l^s+q_k)}\right)\\
\nonumber\\
\mu_k:=\sqrt{1+[\max\left(q_k,|\delta_k(q_k+1-\l^s)-q_k|\right)]^2}
\ear
and for any $r > 0$.\\
The assumption $v$), follows by the linear dependence on the state of the system of
(\ref{delta-jacob}).
We can construct an atlas over $\M_{\delta,k}$ by taking the parametrization (\ref{delta-chart}) 
for the points which are not vertices and by defining compatible charts on the vertices so that
$$
Lip_{M_{\delta,k}}\leq 1.
$$
Therefore, Theorem \ref{TECH-REG-2} implies that $U_{\delta,k}(\cdot,\cdot)$ is Lipschitz continuous over $X_k$ and 
its Lipschitz constant is estimated by $\l_k\mu_k$. Now observe that
$$
U_k(S,B)=\bigcup_\delta U_{\delta,k}(S,B),
$$
for every $(S,B)\in X_k$. Therefore by the definition of the metric $d_H$, we have
$$
U_{\delta,k}(S,B)\subset \left(U_{\delta,k}(S^\p,B^\p)\right)^{\s}\subset \left(U_{k}(S^\p,B^\p)\right)^\s,
$$
which yields
$$
U_k(S,B) \subset\left(U_{k}(S^\p,B^\p)\right)^\s,
$$
for every $(S,B),(S^\p,B^\p)\in X_k$ and $\s>\l_k\mu_k|(S-S^\p,B-B^\p)|$. This implies the Lipschitz regularity
of $U_k(\cdot,\cdot)$. \\
To prove the Lipschitz regularity of $J_k$, it suffices proving that assumption $2)$ in Theorem
\ref{LIP} holds true. In fact we can repeat the proof of Theorem \ref{LIP} by replacing assumption
$1)$ by the Lipschitz regularity of $U_k(\cdot,\cdot)$.\\
By (\ref{model.1}), we have
\bar
f_k(S,B,u,v,y^s,y^b)=\left([S(1-u)+Bv(1-\l^s)]y^s,[B(1-v)+Su(1-\l^b)]y^b\right),
\ear
for each $(S,B,u,v,y^s,y^b)\in X_k\times\M_k\times \Y_k$. Using the linearity and (\ref{xk}), it can be proved that 
$f_k(\cdot,\cdot,\cdot,\cdot,y^s,y^b)$ is Lipschitz continuous over $X_k\times\M_k$ and its
Lipschitz constant is given by
$$
\sqrt{1+D_k^2+\left(1-\min(\l^s,\l^b)\right)^2}\max(y^s,y^b),
$$
which is $p_k$-integrable.
\cvd
\end{Dim}

% --------------------------------------------------------
\section{Acknowledgments}
% --------------------------------------------------------
The authors wish to thank Prof. Wolfgang Runggaldier for his support and his comments, and
Prof. Benedetto Piccoli for his encouragement and his patience in reading countless versions of
the work.\\
The research was supported by the Italian National Research Council and the work of the first author
was carried out while he was visiting the Department of Pure and Applied Mathematics at University
of Padova (Italy).

% --------------------------------------------------------
% BIBLIOGRAPHY
% --------------------------------------------------------

\end{document}